\documentclass[11pt, a4paper]{amsart}
\usepackage{amsfonts}

\usepackage[utf8]{inputenc}
\usepackage{diffcoeff}
\usepackage{xcolor}
\newtheorem{theorem}{\bf Theorem}[section]
\newtheorem{remark}{\bf Remark}[section]
\newtheorem{proposition}{Proposition}[section]
\newtheorem{lemma}{Lemma}[section]
\newtheorem{corollary}{Corollary}[section]

\usepackage{enumerate,mathrsfs}
\usepackage{amssymb,amsmath,graphicx,amsthm,mathtools,hyperref}

\usepackage{booktabs}

\theoremstyle{plain}
\usepackage{siunitx}

\usepackage{algorithm}
\usepackage{algorithmic,colortbl}

\makeindex
\usepackage[intoc]{nomencl} 

\makenomenclature

\usepackage{hyperref}
\hypersetup{nesting=true,debug=true,naturalnames=true}
\usepackage{graphicx,amssymb,upref}

\usepackage{subcaption}
\usepackage{enumerate,float}
\usepackage{fullpage}
\makeatletter
\@namedef{subjclassname@2020}{\textup{2020} Mathematics Subject Classification}
\makeatother

\usepackage{lipsum}
\begin{document}

\title{Multivariate Tempered Space-Fractional Negative binomial Process and Risk Models with Shocks}
	\author[Ashok Kumar Pathak]{Ashok Kumar Pathak}
\author{Ritik Soni$^{*}$\\\\ 
	D\lowercase{epartment of} M\lowercase{athematics and} S\lowercase{tatistics}, C\lowercase{entral} U\lowercase{niversity of} P\lowercase{unjab},\\
	B\lowercase{athinda}, P\lowercase{unjab}-151401, I\lowercase{ndia}.\\
E\lowercase{-mail} A\lowercase{ddress:  ashokiitb09@gmail.com, ritiksoni2012@gmail.com*}}

\thanks{The research of  R. Soni was supported by CSIR, Government of India.}


			\begin{abstract} In this paper, we first define the multivariate tempered space-fractional Poisson process (MTSFPP) by time-changing the multivariate Poisson process with an independent tempered $\alpha$-stable subordinator. Its distributional properties, the mixture tempered time and space variants  and their PDEs connections are studied. Then we define the multivariate tempered space-fractional negative binomial process (MTSFNBP) and explore its key features. The L\'{e}vy measure density for the MTSFNBP is also derived. We present a bivariate risk model with a common shock driven by the tempered space-fractional negative binomial process. We demonstrate that the total claim amount process is stochastically equivalent to a univariate generalized Cramer-Lundberg risk model. In addition, some important ruin measures such as ruin probability, joint distribution of time to ruin and deficit at ruin along with governing integro-differential equations are obtained. Finally we show that the underlying risk process exhibits the long range dependence property.			
\end{abstract}
			\maketitle
			
\section{Introduction}
The Poisson process is the most commonly used point process for modeling counting phenomena.  It has several limitations due to its light-tailed waiting time. Addressing these limitations, researchers have emphasized on the fractional generalizations of the Poisson process via subordination techniques which have several applications in diverse disciplines such as insurance, economics, reliability theory, medical, queuing theory, hydrology and so on (see \cite{Beghin2018, Biard2014, Guler2022, Kumar2020, Laskin2003}). The subordinated variants of the Poisson process are mainly classified into two categories,
namely, the time-fractional and the space-fractional Poisson processes. The time-fractional
Poisson process (TFPP) is obtained from the homogeneous Poisson process (HPP)  by time changing it with an independent inverse $\alpha$-stable subordinator, whereas the space-fractional Poisson process (SFPP) is the subordinated form of the HPP with an independent  $\alpha$-stable subordinator (see   Meerschaert \textit{et al.} \cite{Meerschaert2011} and Orsingher and Polito \cite{Orsingher}). For more recent development on subordinated processes, one can refer to \cite{Gupta2023, Leonenko2015, Maheshwari2019} and \cite{SoniSPL2024}. Recently, Vellaisamy and Maheshwari \cite{Vellaisamy2018} considered a fractional negative binomial process (FNBP) by time changing TFPP with respect to an independent  gamma subordinator. The tempered and space-fractional variants of the FNBP have also been studied in the literature (see \cite{Beghin2018, Maheshwari}). By tempering a stable process, large jumps become significantly smaller, whereas small jumps retain their original stable-like behavior, resulting in integer order moments (see \cite {Maheshwari2023} and \cite{Rosnski2007}).

The multivariate extension of the counting processes have received considerable attention by researchers in recent years for modeling various complex random phenomena that arise  in our daily life. Despite that, there is limited research  on the multivariate generalization of the fractional Poisson process although interest in this topic has risen significantly. By subordination, one can introduce several multivariate counting processes which can be obtained  by time-changing a multivariate L\'{e}vy process with a one-dimensional increasing L\'{e}vy process or through multivariate subordination (see \cite{Barndorff2001} and \cite{Sato}). Next we briefly outline the key studies conducted in this area. Beghin and Macci \cite{Beghin2016} discussed a multivariate fractional Poisson process by considering common random time-change of a finite-dimensional independent Poisson process. Beghin and Macci \cite{Beghin2017} obtained  results related to large deviations for a multivariate version of the alternating fractional Poisson process.  Beghin and Vellaisamy \cite{Beghin2018} investigated the time-changed variant of the  multidimensional space-fractional Poisson process using a common independent gamma subordinator. Recently, Di Crescenzo and Meoli \cite{DiCrescenzo2023} studied competing risks and shock models governed by a bivariate space-fractional Poisson process and explored its ageing notions. Ritik \textit{et al.} \cite{SoniBTFPP} extended the work of Di Crescenzo and Meoli \cite{DiCrescenzo2023} by employing the tempered $\alpha$-stable subordinator and discussed a generalized shock model as well as reliability related findings.

In this paper, we first introduce a multivariate tempered space-fractional Poisson process (MTSFPP) by time-changing the multivariate Poisson process with an independent tempered $\alpha$-stable subordinator and study its important properties. Additionally, we define the mixture tempered time and mixture tempered space variants of the  MTSFPP and their connections to PDEs are also explored. We introduce the time changed variant of MTSFPP with respect to an independent gamma subordinator and call the resulting process  as multivariate tempered space-fractional negative binomial process (MTSFNBP). This can be viewed as a multivariate generalization of the tempered space-fractional negative binomial process (TSFNBP) and a tempered version of the multivariate space fractional negative binomial process studied by Maheshwari \cite{Maheshwari2023} and Beghin and Vellaisamy \cite{Beghin2018}, respectively. We derive its probability mass function (pmf), probability generating function (pgf) and  their governing
differential equations. The L\'{e}vy measure density for the MTSFNBP is also obtained. We propose a bivariate risk model with common shocks based on the tempered space-fractional negative binomial process as an application in risk theory.  It is shown that the total
claim amount process is stochastically equivalent to a univariate risk model governed by the compound analogue of the TSFNBP. This makes it possible to reduce the bivariate model under consideration to a generalized Cramer-Lundberg risk model. Some important ruin measures such as ruin probability, joint distribution of time to ruin and deficit at ruin, and governing integro-differential equations are derived. It has been shown finally that the underlying risk process exhibits the long range dependence (LRD) property.

The structure of the paper is as follows: In Section 2, we present some preliminary notations
and definitions. In Section 3, we introduce the multivariate tempered space-fractional Poisson
process (BTSFPP) and discuss its mixture tempered versions. The multivariate tempered space-fractional negative binomial process and its main characteristics are studied in Section 4.  As an application in ruin theory, a bivariate risk model 
governed by tempered space-fractional negative binomial process is proposed and some results related to ruin theory are presented in Section 5. Finally, some concluding remarks are discussed in the last section.
\section{Preliminaries} 
		In this section, we provide some notations, definitions and results which are important for the development of the results in the subsequent sections.
Let $\mathbb{N}$ denotes the set of natural numbers and $\mathbb{N}_0 = \mathbb{N} \cup \{0\}$. Let $\mathbb{R}$ and $\mathbb{C}$ denote the set of real numbers and set of complex numbers, respectively. 

\subsection{Special Functions and Derivatives}\hfill \\
(i) The three-parameters Mittag-Leffler function $E_{\beta, \gamma}^{\alpha}(z)$ is defined as (see \cite{Podlubny})
\begin{equation*}\label{ml123}
E_{\beta, \gamma}^{\alpha}(z) = \sum_{k=0}^{\infty} \frac{z^k}{k!\Gamma(\gamma+\beta k)}\frac{\Gamma(\alpha +k)}{\Gamma(\alpha)}, \;\;\beta, \gamma, \alpha,  z \in \mathbb{C} \text{ and } \text{Re}(\beta)>0, \text{Re}(\gamma)>0, \text{Re}(\alpha)>0.
\end{equation*}
(ii) The generalized Wright function is defined by (see \cite{Kilbas})
\begin{equation}\label{gwf11}
_p\psi_q \left[z\; \vline \;\begin{matrix}
\left(\alpha_i, \beta_i\right)_{1,p}\\
(a_j,b_j)_{1,q}
\end{matrix} \right] = \sum_{k=0}^{\infty} \frac{z^k}{k!} \frac{\prod_{i=1}^{p} \Gamma(\alpha_i + \beta_i k)}{\prod_{j=1}^{q}\Gamma(a_j + b_j k)},\;\;  z, \alpha_i, a_i \in \mathbb{C}\; \text{and}\; \beta_i, b_i \in \mathbb{R},
\end{equation}
under the convergence condition \begin{equation*}
\sum_{j=1}^{q} b_j - \sum_{i=1}^{p} \beta_i >-1.
\end{equation*}
(iii) Let $f :[a,b] \subset \mathbb{R}\longrightarrow\mathbb{R}$ be $(n+1)$-times continuous differentiable function for $n < \tau <n+1$. Then, the Riemann-Liouville fractional derivative of order $\tau>0$ is defined as (see \cite{Podlubny})
\begin{equation*}
_aD^\tau_tf(t)=\bigg( \frac{d}{dt}\bigg)^{n+1}\int_{a}^{t}(t-u)^{n-\tau}f(u)du.
\end{equation*}

\subsection{L\'{e}vy Subordinator} A L\'{e}vy subordinator denoted by $\{S (t)\}_{t \geq 0}$ is a  non-decreasing L\'{e}vy process with Laplace transform (see \cite[Section 1.3.2]{Applebaum})
\begin{equation*}
\mathbb{E}\left(e^{-uS(t)}\right) = e^{-t \psi(u)},\;\; u \geq 0,
\end{equation*}
where $\psi(u)$ is Laplace exponent given by (see \cite[Theorem 3.2]{Schilling})
\begin{equation*}
\psi(u) = \eta u +\int_{0}^{\infty} (1-e^{-ux}) \nu (dx), \;\; \eta \geq 0.
\end{equation*}
Here $\eta$ is the drift coefficient and $\nu$ is a non-negative L\'{e}vy measure on the positive half-line satisfying
\begin{equation*}
\int_{0}^{\infty} \min \{x,1\} \nu (dx) < \infty, \;\;\; \text{ and } \;\;\; \nu ([0,\infty)) =\infty,
\end{equation*}
so that $\{S (t)\}_{t \geq 0}$ has strictly increasing sample paths almost surely (a.s.) (for more details see  \cite[Theorem 21.3]{Sato}).
\subsubsection{\textbf{Tempered $\alpha$-Stable Subordinator}} For $\alpha \in (0,1)$ and $\theta >0$, the tempered $\alpha$-stable subordinator $\{S^{\alpha, \theta}(t)\}_{t\geq 0}$ is defined by the Laplace transform (see \cite{Kumar})
\begin{equation*}\label{tss11}
\mathbb{E}[e^{-u S^{\alpha, \theta}(t) }] =  e^{\displaystyle -t\left((u+\theta)^\alpha - \theta^\alpha \right)},
\end{equation*} 
with Laplace exponent $\psi(u) = (u+\theta)^\alpha - \theta^\alpha$.
\subsubsection{\textbf{Gamma Subordinator}} Let $\Gamma(t)\sim G(\mu, \rho t), \mu >0, \rho > 0$ be the gamma subordinator with probability density function (pdf) given by
\begin{equation}\label{gamma1}
f_G(x,t) = \frac{\mu^{\rho t}}{\Gamma(\rho t)}x^{\rho t-1}e^{-\mu x}, \;\; x>0.
\end{equation}
Its Laplace transform takes the following form
\begin{equation}\label{gammalt}
\mathbb{E}[e^{-u\Gamma(t)}] = \left(1+\frac{u}{\mu}\right)^{-\rho t}, \;\; u>0.
\end{equation}

\noindent The density in (\ref{gamma1}) of the gamma subordinator $\{\Gamma(t)\}_{t \geq 0}$ satisfies the cauchy problem (see \cite{Beghin2014})
\begin{equation*}
\diffp{}{x} f_G(x,t) = -\mu \left(1-e^{-\frac{1}{\rho}\partial_t}\right) f_G(x,t),
\end{equation*}
with $f_G(x,0) = \delta(x)$ and $\lim_{|x| \rightarrow \infty}f_G(x,t) =0.$ Here, $\delta(x)$ denotes the Dirac delta function and $e^{-\frac{1}{\rho}\partial_t}$ is the partial version of the shift operator defined as
\begin{equation*}
e^{-\frac{1}{\rho}\partial_t}f(t) = \sum_{n=0}^{\infty} \left(-\frac{1}{\rho}\right)^n\frac{\partial_t^n}{n!} f(t) = f(t-1/\rho),
\end{equation*}
where $\partial_t = \diffp{}{t}$ is defined for any analytic function on $\mathbb{R}.$ 
\subsection{Tempered Space-Fractional Poisson Process}
Let $\{\mathcal{N}(t,\lambda)\}_{t \geq 0}$ be the homogeneous Poisson process with parameter $\lambda >0$. Then, the tempered space-fractional Poisson process (TSFPP) denoted by $\{\mathcal{N}^{\alpha, \theta}(t,\lambda)\}_{t \geq 0}$ is defined by time-changing the HPP with an independent tempered $\alpha$-stable subordinator (TSS) as (see \cite{Gupta})
\begin{equation*}
\mathcal{N}^{\alpha, \theta}(t,\lambda) := \mathcal{N}(S^{\alpha, \theta}(t), \lambda).
\end{equation*}
Its pmf $p^{\alpha, \theta}(k,t)$ is given by (see \cite[Formula (26)]{Gupta1})
\begin{equation*}
p^{\alpha, \theta}(k,t) = \frac{(-1)^k}{k!}  e^{t\theta^{\alpha}}\sum_{i=0}^{\infty} \frac{\theta^i}{\lambda^i i!} \; _1\psi_1 \left[-\lambda^\alpha t\; \vline \;\begin{matrix}
\left(1,\alpha\right)\\
(1-k-i, \alpha)
\end{matrix} \right].
\end{equation*}
\subsection{Tempered Space-Fractional Negative Binomial Process}
Let $\{\mathcal{N}^{\alpha, \theta}(t,\lambda)\}_{t \geq 0}$ is TSFPP and $\Gamma(t)\sim G(\mu, \rho t)$ be a gamma subordinator. Then, the tempered space-fractional negative binomial process (TSFNBP) is defined as (see \cite{Maheshwari2023})
\begin{equation*}
\mathcal{Q}(t) := \mathcal{N}^{\alpha, \theta}(\Gamma(t),\lambda).
\end{equation*}
 Its pgf and mean are of the following form (see \cite[Section 4.2]{Maheshwari2023})
 \begin{equation*}
 \mathbb{E}[u^{\mathcal{Q}(t)}] = \left(1+\mu^{-1}\left[\left(\lambda (1-u)+\theta\right)^{\alpha}-\theta^\alpha\right]\right)^{-\rho t},
 \end{equation*}
 with $\mathbb{E}[\mathcal{Q}(t)] = \alpha \lambda \rho t \mu^{-1} \theta^{\alpha-1}$.

\subsection{Backward Shift Operators}
Let $B$ be the backward shift operator defined by $B[\xi(k)] = \xi(k-1)$. For the fractional difference operator $(I-B)^\alpha$, we have  (see \cite{Orsingher} and \cite[page 91]{Tsay})
\begin{equation*}\label{do1}
(I-B)^\alpha = \sum_{i=0}^{\infty} \binom{\alpha}{j}(-1)^i B^i, \;\;\alpha \in (0,1),
\end{equation*}
where $I$ is an identity operator.\\
Furthermore, let  $\{B_i\}$, $ i\in \{1,2,\dots, m\}$ be the operators defined as 
\begin{equation*}\label{do2}
B_i[\xi(k_1, k_2, \dots, k_m)] = \xi(k_1, k_2, \dots, k_{i}-1,\dots, k_m).
\end{equation*}
When $m=1$,  $B_i$'s act same as the operator $B$.  

\subsection{Long-Range Dependence} 
	For $0 <s <t,$ let the correlation function Corr$[X(s), X(t)]$ for a non-stationary stochastic process $\{X(t)\}_{t \geq 0}$ satisfes (see \cite[Definition 2]{Kumar2020} and \cite{Nane})
\begin{equation*}
\text{ Corr}[X(s), X(t)] \sim c(s)t^{-d},\;\; \text{as} \; t \rightarrow \infty,
\end{equation*}
for some $c(s) >0$ and $d > 0$.
The process $\{X(t)\}_{t \geq 0}$ is said to have the LRD property if $d \in (0,1)$.

\section{Multivariate Tempered Space-Fractional Poisson Process}
Let $\mathcal{N}(t) = (N_1(t, \lambda_1), N_2(t, \lambda_2), \dots, N_m(t, \lambda_m)),\; t \geq 0 $ be a 
multivariate HPP such that $\{N_i(t, \lambda_i)\}_{t \geq 0}$, $i=1,2,\ldots,m$ are independent HPPs with  parameters $\lambda_i$. We define the multivariate tempered space-fractional Poisson process (MTSFPP) by time-changing the multivariate HPP with an independent TSS, $\{S^{\alpha, \theta}(t)\}_{t \geq 0}$ as 
\begin{equation}\label{bp1}
\mathcal{N}^{\alpha, \theta}(t) :=\left( \mathcal{N}_1^{\alpha, \theta}(t,\lambda_1), \mathcal{N}_2^{\alpha, \theta}(t,\lambda_2),\ldots, \mathcal{N}_m^{\alpha, \theta}(t,\lambda_m)\right), \;\; t \geq 0,
\end{equation}
where $\mathcal{N}_i^{\alpha, \theta}(t,\lambda_i) = N_i\left(S^{\alpha, \theta}(t), \lambda_i\right), \; i=1,2,\ldots,m$ are i.i.d. TSFPPs.

We denote any arbitrary multivariate vector of constants by $\textbf{a}=(a_1,\ldots, a_m)$, where  $a_1,\ldots, a_m$ are nonnegative integers.
Let  $\textbf{b}=(b_1,\ldots, b_m)$  and $\textbf{0} = (0,\ldots,0)$ be the null vector.  We write $\textbf{a} \geq  \textbf{b}$ (or $\textbf{a} \leq  \textbf{b}$)  to mean that $a_i\geq b_i$ (or $a_i\leq b_i$) for $i=1,\ldots,m$. Further, we denote $\textbf{k}=(k_1, \ldots,k_m)$, $\textbf{r}=(r_1,\ldots, r_m)$ and $ \boldsymbol{\lambda}=(\lambda_1, \ldots,\lambda_m)$. Finally, we consider $\textbf{u}=(u_1,\ldots, u_m)\in [0, 1]^m$ for the pgf existence.

\begin{proposition}
	For $\textbf{u} = (u_1,u_2,\ldots, u_m) \in [0,1]^m$, the pgf  $G^{\alpha, \theta}(\textbf{u};t)$ for the MTSFPP is given by
	\begin{equation*}
	G^{\alpha, \theta}(\textbf{u};t) = \exp\left\{\displaystyle
		-t\left(\left[\sum_{i=1}^{m}\lambda_i(1-u_i) + \theta\right]^{\alpha} - \theta^\alpha\right)\right\},
	\end{equation*}
	and it satisfies the differential equation of the following form
	\begin{equation}\label{de23}
	\diffp{}{t}	G^{\alpha, \theta}(\textbf{u};t) = -\left(\left[\sum_{i=1}^{m}\lambda_i(1-u_i) + \theta\right]^{\alpha} - \theta^\alpha\right)	G^{\alpha, \theta}(\textbf{u};t), \;\; 	G^{\alpha, \theta}(\textbf{u};0) =1.
	\end{equation}
\end{proposition}
\begin{proof}
	For $\lambda >0$, the pgf for the TSFPP is derived as (see \cite{Gupta1})
	\begin{align*}
	\mathbb{E}\left[u^{ \mathcal{N}^{\alpha, \theta}(t,\lambda)}\right] 
	&= \mathbb{E}\left[\mathbb{E}[u^{\mathcal{N}(S^{\alpha, \theta}(t), \lambda)}|S^{\alpha, \theta}(t)]\right]\\
	&= \mathbb{E}\left[e^{-\lambda (1-u)S^{\alpha, \theta}(t)}\right]\\
	&= e^{\displaystyle -t((\lambda(1-u)+\theta)^\alpha - \theta^\alpha)}.
	\end{align*}
	We define the pgf of MTSFPP as
	\begin{equation*}
	G^{\alpha, \theta}(\textbf{u};t) =	\mathbb{E}\left[\textbf{u}^{\mathcal{N}^{\alpha, \theta}(t)}\right] = \sum_{\textbf{k} \geq \textbf{0}}^{} u_1^{k_1} u_2^{k_2}\cdots u_m^{k_m} 	q^{\alpha, \theta}(\textbf{k},t),
	\end{equation*}
	where $q^{\alpha, \theta}(\textbf{k},t)$ denotes the pmf of the MTSFPP.\\
	By conditioning argument, we get
	\begin{align*}\label{pgf1123}
	G^{\alpha, \theta}(u;t) =& \mathbb{E}\left[\mathbb{E}[u^{\mathcal{Q}^{\alpha, \theta}(t, \lambda)}\;|\;S^{\alpha, \theta}(t)]\right]\nonumber\\=& \mathbb{E}\left[e^{\sum_{i=1}^{m}\lambda_i(u_i-1)S^{\alpha, \theta}(t)}\right] \nonumber\\=& \exp\left\{\displaystyle
	-t\left(\left[\sum_{i=1}^{m}\lambda_i(1-u_i) + \theta\right]^{\alpha} - \theta^\alpha\right)\right\}.
	\end{align*}
	By calculus, we obtain (\ref{de23}) and   the condition trivially holds for $t=0$.
\end{proof}
The pmf of the MTSFPP and associated differential equations are obtained in the following results.
\begin{proposition}
	For $\alpha \in (0,1)$ and $\textbf{k} \geq \textbf{0}$, the pmf $q^{\alpha, \theta}(\textbf{k},t) =Pr\{\mathcal{N}^{\alpha, \theta}(t) =\textbf{k}\}$ is given by
	\begin{equation}\label{pmf2}
	q^{\alpha, \theta}(\textbf{k},t) = \left(-\frac{1}{S({\boldsymbol{\lambda}})}\right)^{S(\boldsymbol{k})}\prod_{j=1}^{m} \frac{\lambda_j^{k_j}}{k_j!}e^{t \theta ^\alpha}\sum_{i=0}^{\infty} \frac{\theta^i}{i! (S(\boldsymbol{\lambda}))^i} \; _1\psi_1 \left[-(S(\boldsymbol{\lambda}))^\alpha t\; \vline \;\begin{matrix}
	\left(1,\alpha\right)\\
	(1-S(\boldsymbol{k})-i, \alpha)
	\end{matrix} \right],
	\end{equation}
	satisfies the following differential equation
		\begin{equation}\label{de111}
		\frac{d}{dt} q^{\alpha, \theta}(\textbf{k},t) = -(S(\boldsymbol{\lambda}))^\alpha \left( \left(I- \frac{\sum_{i=1}^{m}\lambda_i B_i - \theta}{S(\boldsymbol{\lambda})}\right)^\alpha -\left( \frac{\theta}{S(\boldsymbol{\lambda})}\right)^\alpha\right)q^{\alpha, \theta}(\textbf{k},t), \;\;  q^{\alpha, \theta}(\textbf{0},t) = 1.
		\end{equation} 
\end{proposition}
\begin{proof}
	The proof of the proposition follow on the same path to that of the Proposition 3.1 and Theorem 3.2 in \cite{SoniBTFPP}. The convergence of $\;_1\Psi_1$ holds due to condition in (\ref{gwf11}).
\end{proof}

\begin{remark}
	When $\theta=0$ and $\alpha=1$, (\ref{de111}) reduces to
	\begin{equation}\label{mhpp1}
		\frac{d}{dt} q(\textbf{k},t) = -S(\boldsymbol{\lambda})  \left(I- \frac{\sum_{i=1}^{m}\lambda_i B_i}{S(\boldsymbol{\lambda})}\right) q(\textbf{k},t), \;\;  q(\textbf{0},t) = 1,
	\end{equation}
	which coincides with the pmf $q(\textbf{k},t) = Pr\{\mathcal{N}(s) =\textbf{k}\}$ of the multivariate HPP reported in \cite{Beghin2016}.
\end{remark}
\begin{theorem}
	The pmf $q^{\alpha, \theta}(\textbf{k},t)$ satisfies the following differential equation
	\begin{equation*}
		\left(\theta^\alpha - \frac{d}{dt}\right)^{1/\alpha}q^{\alpha, \theta}(\textbf{k},t) = \theta q^{\alpha, \theta}(\textbf{k},t) + S(\boldsymbol{\lambda})  \left(I- \frac{\sum_{i=1}^{m}\lambda_i B_i }{S(\boldsymbol{\lambda})}\right) q^{\alpha, \theta}(\textbf{k},t).
	\end{equation*}
\end{theorem}
\begin{proof}
	Let $h^{\alpha, \theta}(x,t)$ be the pdf of the TSS satisfies (see \cite{Beghin2015})
	\begin{equation*}
	\diffp{}{x}h^{\alpha, \theta}(x,t) = -\theta h^{\alpha, \theta}(x,t) + \left(\theta^\alpha - \frac{d}{dt}\right)^{1/\alpha}h^{\alpha, \theta}(x,t), \;\; h^{\alpha, \theta}(x,0) = \delta(x).
	\end{equation*}
	With the help of (\ref{mhpp1}), we get
	\begin{align*}
		\left(\theta^\alpha - \frac{d}{dt}\right)^{1/\alpha}q^{\alpha, \theta}(\textbf{k},t) &= \int_{0}^{\infty}  Pr\{\mathcal{N}(s) =\textbf{k}\} \left(\theta^\alpha - \frac{d}{dt}\right)^{1/\alpha} h^{\alpha, \theta}(s,t)ds\\
		&= \theta q^{\alpha, \theta}(\textbf{k},t) - \int_{0}^{\infty} \frac{d}{ds}Pr\{\mathcal{N}(s) =\textbf{k}\}  h^{\alpha, \theta}(s,t) ds\\
		&=  \theta q^{\alpha, \theta}(\textbf{k},t) + \int_{0}^{\infty} S(\boldsymbol{\lambda})  \left(I- \frac{\sum_{i=1}^{m}\lambda_i B_i }{S(\boldsymbol{\lambda})}\right)Pr\{\mathcal{N}(s) =\textbf{k}\}  h^{\alpha, \theta}(s,t) ds\\
		&= \theta q^{\alpha, \theta}(\textbf{k},t) + S(\boldsymbol{\lambda})  \left(I- \frac{\sum_{i=1}^{m}\lambda_i B_i }{S(\boldsymbol{\lambda})}\right) q^{\alpha, \theta}(\textbf{k},t).
	\end{align*}
	Hence, the theorem is proved.
\end{proof}
\subsection{Multivariate Mixture Tempered Space-Fractional Poisson Process} \hfill \\
For $\alpha_1, \alpha_2 \in (0,1)$ and $\theta_1, \theta_2 >0$. Let $\{\mathcal{S}_{\alpha_1, \theta_1}^{\alpha_2, \theta_2}(t)\}_{t \geq 0}$ be the mixture of tempered stable subordinator (MTSS). Then, we define the multivariate mixture tempered space-fractional Poisson process (MMTSFPP) $\{M_1(t)\}_{t \geq 0}$ as
\begin{equation}
M_1(t) = \left(N_1(\mathcal{S}_{\alpha_1, \theta_1}^{\alpha_2, \theta_2}(t), \lambda_1), N_2(\mathcal{S}_{\alpha_1, \theta_1}^{\alpha_2, \theta_2}(t), \lambda_2), \dots, N_m(\mathcal{S}_{\alpha_1, \theta_1}^{\alpha_2, \theta_2}(t), \lambda_m)\right), \;\; t \geq 0,
\end{equation}
where $\{N_i(t, \lambda_i)\}_{t \geq 0}, \; i=1,2,\ldots, m$ are HPPs independent from $\{\mathcal{S}_{\alpha_1, \theta_1}^{\alpha_2, \theta_2}(t)\}_{t \geq 0}$.\\
Now, we derive the governing differential equation to the pmf $\hat{p}_1(\textbf{k},t)$ of the MMTSFPP.
\begin{proposition}
	The pmf $\hat{p}_1(\textbf{k},t)$ satisfies the following differential equation
\begin{align*}
	\frac{d}{dt} \hat{p}_1(\textbf{k},t) &=  \int_{0}^{\infty} Pr\{\mathcal{N}(s) =\textbf{k}\} \left(	-	\eta_1\left(\theta_1 +\diffp{}{s}\right)^{\alpha_1}h_{\alpha_1, \theta_1}^{\alpha_2, \theta_2}(s,t) -\eta_2\left(\theta_2 +\diffp{}{s}\right)^{\alpha_2}h_{\alpha_1, \theta_1}^{\alpha_2, \theta_2}(s,t)\right)ds\\
		&\;\;\; +  \eta_1 \hat{p}_1(\textbf{k},t) +\theta_2^{\alpha_1} \eta_2 \hat{p}_1(\textbf{k},t).
\end{align*}
\end{proposition}

\begin{proof}
		The pdf $h_{\alpha_1, \theta_1}^{\alpha_2, \theta_2}(x,t) $ of the MTSS satisfies the following differential equation (see \cite{Gupta2021})
		\begin{align*}
		\diffp{}{t} h_{\alpha_1, \theta_1}^{\alpha_2, \theta_2}(x,t) &+ \eta_1\left(\theta_1 +\diffp{}{x}\right)^{\alpha_1}h_{\alpha_1, \theta_1}^{\alpha_2, \theta_2}(x,t) + \eta_2\left(\theta_2 +\diffp{}{x}\right)^{\alpha_2}h_{\alpha_1, \theta_1}^{\alpha_2, \theta_2}(x,t)\\
		&= \theta_1^{\alpha_1} \eta_1 h_{\alpha_1, \theta_1}^{\alpha_2, \theta_2}(x,t) + \theta_2^{\alpha_1} \eta_2 h_{\alpha_1, \theta_1}^{\alpha_2, \theta_2}(x,t).
		\end{align*}
		By using (\ref{mhpp1}), we have
		\begin{align*}
		\frac{d}{dt} \hat{p}_1(\textbf{k},t) 
			&= \int_{0}^{\infty} Pr\{\mathcal{N}(s) =\textbf{k}\} 	\diffp{}{t}h_{\alpha_1, \theta_1}^{\alpha_2, \theta_2}(s,t)ds\\
			&=  \int_{0}^{\infty} Pr\{\mathcal{N}(s) =\textbf{k}\} \left(	-	\eta_1\left(\theta_1 +\diffp{}{s}\right)^{\alpha_1}h_{\alpha_1, \theta_1}^{\alpha_2, \theta_2}(s,t) -\eta_2\left(\theta_2 +\diffp{}{s}\right)^{\alpha_2}h_{\alpha_1, \theta_1}^{\alpha_2, \theta_2}(s,t)\right)ds\\
			&\;\;\; + \int_{0}^{\infty} Pr\{\mathcal{N}(s) =\textbf{k}\} \theta_1^{\alpha_1} \eta_1 h_{\alpha_1, \theta_1}^{\alpha_2, \theta_2}(s,t)ds + \int_{0}^{\infty} Pr\{\mathcal{N}(s) =\textbf{k}\}\theta_2^{\alpha_1} \eta_2 h_{\alpha_1, \theta_1}^{\alpha_2, \theta_2}(s,t)ds.
		\end{align*}
		On simplifying, we get the proposition.
\end{proof}
\subsection{Multivariate Mixture Tempered Time-Fractional Poisson Process} \hfill \\

For $\alpha_1, \alpha_2 \in (0,1)$ and $\theta_1, \theta_2 >0$. Let $\{\mathcal{E}_{\alpha_1, \theta_1}^{\alpha_2, \theta_2}(t)\}_{t \geq 0}$ be the inverse of the mixture of tempered stable subordinator (IMTSS). Then, we define the multivariate mixture tempered space-fractional Poisson process (MMTTFPP) $\{M_2(t)\}_{t \geq 0}$ as
\begin{equation}
M_2(t) = \left(N_1(\mathcal{E}_{\alpha_1, \theta_1}^{\alpha_2, \theta_2}(t), \lambda_1), N_2(\mathcal{E}_{\alpha_1, \theta_1}^{\alpha_2, \theta_2}(t), \lambda_2), \dots, N_m(\mathcal{E}_{\alpha_1, \theta_1}^{\alpha_2, \theta_2}(t), \lambda_m)\right), \;\; t \geq 0,
\end{equation}
where $\{N_i(t, \lambda_i)\}_{t \geq 0}, \; i=1,2,\ldots, m$ are HPPs independent from $\{\mathcal{E}_{\alpha_1, \theta_1}^{\alpha_2, \theta_2}(t)\}_{t \geq 0}$.\\
Next, we obtain the governing differential equation to the pmf $\hat{p}_2(\textbf{k},t)$ of the MMTTFPP.
\begin{proposition}
	The pmf $\hat{p}_2(\textbf{k},t)$ satisfies the following differential equation
	\begin{align}\label{de1}
		&\left[	\eta_1\left(\theta_1 +\diffp{}{t}\right)^{\alpha_1} +\eta_2\left(\theta_2 +\diffp{}{t}\right)^{\alpha_2}\right]\hat{p}_2(\textbf{k},t) = -S(\boldsymbol{\lambda})  \left(I- \frac{\sum_{i=1}^{m}\lambda_i B_i }{S(\boldsymbol{\lambda})}\right)\hat{p}_2(\textbf{k},t) \nonumber \\&+ \left[\theta_1^{\alpha_1}\eta_1 + \theta_2^{\alpha_1}\eta_2\right]\hat{p}_2(\textbf{k},t)- \eta_1 t^{-\alpha_1}E_{1,1-\alpha_1}^{1-\alpha_1}(-\theta_1 t)\delta(x)
	 - \eta_2 t^{-\alpha_2}E_{1,1-\alpha_2}^{1-\alpha_2}(-\theta_2 t)\delta(x).
	\end{align}
\end{proposition}

\begin{proof}
	The pdf $g_{\alpha_1, \theta_1}^{\alpha_2, \theta_2}(x,t) $ of the IMTSS satisfies the following differential equation (see \cite{Gupta2021})
	\begin{align*}
	\diffp{}{x}g_{\alpha_1, \theta_1}^{\alpha_2, \theta_2}(x,t) = &-\eta_1\left(\theta_1 +\diffp{}{t}\right)^{\alpha_1}g_{\alpha_1, \theta_1}^{\alpha_2, \theta_2}(x,t) -\eta_2\left(\theta_2 +\diffp{}{t}\right)^{\alpha_2}g_{\alpha_1, \theta_1}^{\alpha_2, \theta_2}(x,t) + \theta_1^{\alpha_1}\eta_1g_{\alpha_1, \theta_1}^{\alpha_2, \theta_2}(x,t)\\
	&+  \theta_2^{\alpha_1}\eta_2g_{\alpha_1, \theta_1}^{\alpha_2, \theta_2}(x,t) - \eta_1 t^{-\alpha_1}E_{1,1-\alpha_1}^{1-\alpha_1}(-\theta_1 t)\delta(x) - \eta_2 t^{-\alpha_2}E_{1,1-\alpha_2}^{1-\alpha_2}(-\theta_2 t)\delta(x),
	\end{align*}
	where $g_{\alpha_1, \theta_1}^{\alpha_2, \theta_2}(x,0) = \delta(x)$.\\
With the help of (\ref{mhpp1}), we have\\ 
$\left[	\eta_1\left(\theta_1 +\diffp{}{t}\right)^{\alpha_1} +\eta_2\left(\theta_2 +\diffp{}{t}\right)^{\alpha_2}\right]\hat{p}_2(\textbf{k},t)$
	\begin{align*}
&= \int_{0}^{\infty} Pr\{\mathcal{N}(s) =\textbf{k}\} \left[	\eta_1\left(\theta_1 +\diffp{}{t}\right)^{\alpha_1} +\eta_2\left(\theta_2 +\diffp{}{t}\right)^{\alpha_2}\right]g_{\alpha_1, \theta_1}^{\alpha_2, \theta_2}(s,t)ds\\
&= \int_{0}^{\infty} Pr\{\mathcal{N}(s) =\textbf{k}\} \diffp{}{s} g_{\alpha_1, \theta_1}^{\alpha_2, \theta_2}(s,t)ds + \left[\theta_1^{\alpha_1}\eta_1 + \theta_2^{\alpha_1}\eta_2\right]\hat{p}_2(\textbf{k},t)\\
&\;\;\; - \eta_1 t^{-\alpha_1}E_{1,1-\alpha_1}^{1-\alpha_1}(-\theta_1 t)\delta(x) - \eta_2 t^{-\alpha_2}E_{1,1-\alpha_2}^{1-\alpha_2}(-\theta_2 t)\delta(x)\\
&= -\int_{0}^{\infty}  \left(\frac{d}{ds} Pr\{\mathcal{N}(s) =\textbf{k}\}\right) g_{\alpha_1, \theta_1}^{\alpha_2, \theta_2}(s,t)ds + \left[\theta_1^{\alpha_1}\eta_1 + \theta_2^{\alpha_1}\eta_2\right]\hat{p}_2(\textbf{k},t)\\
&\;\;\; - \eta_1 t^{-\alpha_1}E_{1,1-\alpha_1}^{1-\alpha_1}(-\theta_1 t)\delta(x) - \eta_2 t^{-\alpha_2}E_{1,1-\alpha_2}^{1-\alpha_2}(-\theta_2 t)\delta(x)\\
&= \int_{0}^{\infty}  \left(-S(\boldsymbol{\lambda})  \left(I- \frac{\sum_{i=1}^{m}\lambda_i B_i }{S(\boldsymbol{\lambda})}\right)\right)Pr\{\mathcal{N}(s) =\textbf{k}\} g_{\alpha_1, \theta_1}^{\alpha_2, \theta_2}(s,t)ds + \left[\theta_1^{\alpha_1}\eta_1 + \theta_2^{\alpha_1}\eta_2\right]\hat{p}_2(\textbf{k},t)\\
&\;\;\; - \eta_1 t^{-\alpha_1}E_{1,1-\alpha_1}^{1-\alpha_1}(-\theta_1 t)\delta(x) - \eta_2 t^{-\alpha_2}E_{1,1-\alpha_2}^{1-\alpha_2}(-\theta_2 t)\delta(x).
	\end{align*}
	On suitably arranging the terms, we get the required differential equation.
\end{proof}

\begin{remark}
	For $m=1$, (\ref{de1}) reduces to the differential equation satisfying the pmf of mixture tempered time-fractional Poisson process studied in \cite{Gupta2021}.
\end{remark}

\section{Multivariate Tempered Space-Fractional Negative Binomial Process}
Let $\Gamma(t) \sim G(\mu, \rho t)$ be the gamma subordinator. We define the multivariate tempered space-fractional negative binomial process (MTSFNBP) by time-changing the MTSFPP with an independent gamma subordinator as
\begin{equation*}
\mathcal{Q}^{\alpha, \theta}(m,t) = (\mathcal{Q}_1^{\alpha, \theta}(t), \mathcal{Q}_2^{\alpha, \theta}(t), \ldots, \mathcal{Q}_m^{\alpha, \theta}(t)),
\end{equation*}
where $\mathcal{Q}_i^{\alpha, \theta}(t) =  \mathcal{N}_i^{\alpha, \theta}(\Gamma(t),\lambda_i), \;i=1,2,\ldots,m$ are mutually independent TSFPP time changed by an independent gamma subordinator.\\
With the help of (\ref{gamma1}) and (\ref{pmf2}) and by an use of Fubini-Tonelli's theorem, we derive the pmf $\mathcal{P}^{\alpha, \theta}(\textbf{k},t) = Pr\{\mathcal{Q}^{\alpha, \theta}(m,t) = \textbf{k}\}$ of the MTSFNBP as
	\begin{align*}
\mathcal{P}^{\alpha, \theta}(\textbf{k},t)
&= \int_{0}^{\infty} q^{\alpha, \theta}(\textbf{k},s) f_G(s,t) ds\\
&= \int_{0}^{\infty} \left(-\frac{1}{S({\boldsymbol{\lambda}})}\right)^{S(\boldsymbol{k})}\prod_{j=1}^{m} \frac{\lambda_j^{k_j}}{k_j!}e^{s \theta ^\alpha}\sum_{i=0}^{\infty} \frac{\theta^i}{i! (S(\boldsymbol{\lambda}))^i}\\
&\;\; \times \sum_{l=0}^{\infty} \frac{(- (S({\boldsymbol{\lambda}}))^\alpha s)^l }{l!} \frac{\Gamma(1+\alpha l)}{\Gamma(1-S(\boldsymbol{k})-i+\alpha l)}f_G(s,t) ds\\
&= \left(-\frac{1}{S({\boldsymbol{\lambda}})}\right)^{S(\boldsymbol{k})}\prod_{j=1}^{m} \frac{\lambda_j^{k_j}}{k_j!}\sum_{i=0}^{\infty} \frac{\theta^i}{i! (S(\boldsymbol{\lambda}))^i} \frac{\mu^{\rho t}}{\Gamma(\rho t)}\\
&\;\; \times \sum_{l=0}^{\infty} \frac{(- (S({\boldsymbol{\lambda}}))^\alpha)^l }{l!} \frac{\Gamma(1+\alpha l)}{\Gamma(1-S(\boldsymbol{k})-i+\alpha l)} \int_{0}^{\infty}e^{-s(\mu- \theta^ \alpha)} s^{\rho t -1}ds\\
&= \left(-\frac{1}{S({\boldsymbol{\lambda}})}\right)^{S(\boldsymbol{k})} \left(\frac{\mu}{\mu- \theta^{\alpha}}\right)^{\rho t}\prod_{j=1}^{m} \frac{\lambda_j^{k_j}}{k_j!}\sum_{i=0}^{\infty} \frac{\theta^i}{i! (S(\boldsymbol{\lambda}))^i} \\
&\;\; \times \sum_{l=0}^{\infty} \frac{(- (S({\boldsymbol{\lambda}}))^\alpha)^l }{l!} \frac{\Gamma(1+\alpha l)}{\Gamma(1-S(\boldsymbol{k})-i+\alpha l)} \frac{\Gamma(\rho t +l)}{(\mu- \theta^\alpha)^{l}},  \;\; \mu > \theta^{\alpha}.
\end{align*}
In terms of the generalized Wright function, we get
	\begin{align}\label{mts1}
\mathcal{P}^{\alpha, \theta}(\textbf{k},t) = \frac{1}{\Gamma(\rho t)} &\left(\frac{\mu}{\mu- \theta^{\alpha}}\right)^{\rho t} \left(\frac{-1}{S(\boldsymbol{\lambda})}\right)^{S(\boldsymbol{k})} \prod_{j=1}^m \frac{\lambda_j^{k_j}}{k_j !} \sum_{i=0}^{\infty} \frac{\theta^i}{i! (S(\boldsymbol{\lambda}))^i }\\
& \times \; _2 \psi_1  \left[\frac{-(S(\boldsymbol{\lambda}))^\alpha}{\mu - \theta^{\alpha}} \; \vline \;\begin{matrix}
\left(1,\alpha\right), (\rho t, 1)\\
(1-S(\boldsymbol{k})-i, \alpha)
\end{matrix} \right], \;\; \mu > \theta^{\alpha}.\nonumber
\end{align}

\begin{remark}
	When $\theta =0$, (\ref{mts1}) reduces to
	\begin{equation}\label{st1}
		\mathcal{P}^{\alpha, 0}(\textbf{k},t) = \frac{1}{\Gamma(\rho t)}  \left(\frac{-1}{S(\boldsymbol{\lambda})}\right)^{S(\boldsymbol{k})} \prod_{j=1}^m \frac{\lambda_j^{k_j}}{k_j !}\; _2 \psi_1  \left[\frac{-(S(\boldsymbol{\lambda}))^\alpha}{\mu} \; \vline \;\begin{matrix}
		\left(1,\alpha\right), (\rho t, 1)\\
		(1-S(\boldsymbol{k}), \alpha)
		\end{matrix} \right],
	\end{equation}
	which is the pmf of multivariate space-fractional negative binomial process  studied in \cite{Beghin2018}.\\
	Also, for $m=1$, (\ref{mts1}) and (\ref{st1}) reduces to the tempered space-fractional negative binomial process (see \cite{Maheshwari2023}) and the space-fractional negative binomial process (see \cite{Beghin2018}), respectively.
\end{remark}
\noindent Next, we obtain the differential equations governed by the pmf of the MTSFNBP.
\begin{theorem}
	The pmf $	\mathcal{P}^{\alpha, \theta}(\textbf{k},t)$ satisfies the following differential equation
	\begin{equation}\label{pmf21}
	\mu \left(1-e^{-\frac{1}{\rho}\partial_t}\right) \mathcal{P}^{\alpha, \theta}(\textbf{k},t) = -(S(\boldsymbol{\lambda}))^\alpha \left( \left(I- \frac{\sum_{i=1}^{m}\lambda_i B_i - \theta}{S(\boldsymbol{\lambda})}\right)^\alpha -\left( \frac{\theta}{S(\boldsymbol{\lambda})}\right)^\alpha\right) \mathcal{P}^{\alpha, \theta}(\textbf{k},t).
	\end{equation}
\end{theorem}
\begin{proof}It is known that
	\begin{equation}\label{def1}
	 \mathcal{P}^{\alpha, \theta}(\textbf{k},t) = \int_{0}^{\infty} 	q^{\alpha, \theta}(\textbf{k},s)f_G(s,t) ds.
	\end{equation}
	Therefore, we have
	\begin{align*}
	e^{-\frac{1}{\rho}\partial_t} \mathcal{P}^{\alpha, \theta}(\textbf{k},t) 
	&= \int_{0}^{\infty} 	q^{\alpha, \theta}(\textbf{k},s) e^{-\frac{1}{\rho}\partial_t}f_G(s,t) ds\\
	&= \int_{0}^{\infty} 	q^{\alpha, \theta}(\textbf{k},s) \left[\frac{1}{\mu} \diffp{}{s}f_G(s,t) + f_G(s,t) \right] ds\\
	&= \mathcal{P}^{\alpha, \theta}(\textbf{k},t) -\frac{1}{\mu} \int_{0}^{\infty} \frac{d}{ds} 	q^{\alpha, \theta}(\textbf{k},s) f_G(s,t)  ds.
	\end{align*}
	Substituting (\ref{de111}) and then arranging the terms, the differential equation is obtained.
\end{proof}
 The following lemma is important to obtain the fractional PDE governing the pmf of the MTSFNBP.
\begin{lemma}\label{lemm1}
	(\cite{Vellaisamy2018}) For $\tau \geq 1$, the governing fractional PDE for the gamma subordinator $\{\Gamma(t)\}_{t\geq 0}$ is given by
	\begin{equation*}
	\diffp{^\tau}{t^\tau}	f_G(x,t) = \rho	\diffp{^{\tau-1}}{t^\tau}\left[\log \alpha +\log y - \psi(\rho 
	t)\right]	f_G(x,t), \;\; y >0 \text{ and }
	f_G(x,0) =0,
	\end{equation*}
	where $\psi(x)$ is the digamma function and $	\diffp{^\tau}{t^\tau}(\cdot)$ is the R-L fractional differential operator.
\end{lemma}

\begin{theorem}
	For $\tau \geq 1$, the pmf $\mathcal{P}^{\alpha, \theta}(\textbf{k},t) $ of the MTSFNBP solves the following fractional differential equation
		\begin{equation*}
	\diffp{^\tau}{t^\tau} \mathcal{P}^{\alpha, \theta}(\textbf{k},t) = \rho	\diffp{^{\tau-1}}{t^{\tau-1}}\left[\left(\log \mu - \psi(\rho t) \right) \mathcal{P}^{\alpha, \theta}(\textbf{k},t)+\int_{0}^{\infty} q^{\alpha, \theta}(\textbf{k},s) (\log s) f_G(s,t)ds\right].
	\end{equation*}
\end{theorem}
	\begin{proof}
		Operating the Riemann-Liouville fractional derivative in (\ref{def1}), we get
		\begin{align*}
		\diffp{^\tau}{t^\tau}  \mathcal{P}^{\alpha, \theta}(\textbf{k},t)
		&= \diffp{^\tau}{t^\tau} \int_{0}^{\infty} 	q^{\alpha, \theta}(\textbf{k},s)	f_G(s,t) ds\\
		&= \int_{0}^{\infty} 	q^{\alpha, \theta}(\textbf{k},s) 	\diffp{^\tau}{t^\tau} 	f_G(s,t) ds\\
		&= \int_{0}^{\infty} 	q^{\alpha, \theta}(\textbf{k},s) \left[\rho \diffp{^{\tau-1}}{t^\tau}\left[\log \mu +\log s - \psi(\rho t)\right]	f_G(s,t)\right]ds \;\;\;(\text{by Lemma \ref{lemm1}})\\
		&= \rho	\diffp{^{\tau-1}}{t^{\tau-1}}\int_{0}^{\infty} 	q^{\alpha, \theta}(\textbf{k},s) \left(\log \mu - \psi(\rho t) \right) f_G(s,t)ds\\
		&\;\; +\rho \int_{0}^{\infty}	q^{\alpha, \theta}(\textbf{k},s) (\log s) 	\diffp{^{\tau-1}}{t^{\tau-1}}f_G(s,t)ds.
		\end{align*}
		With the help of simple algebra, the result follows.
	\end{proof}

In the next, we derive the pgf of the MTSFNBP and its governing differential equation.

\begin{theorem}
For $\alpha \in (0,1)$ and $\theta >0$, the pgf $\mathcal{G}^{\alpha, \theta}(\textbf{u},t)$ of the MTSFNBP is given by
\begin{equation}\label{mtsfpgf}
\mathcal{G}^{\alpha, \theta}(\textbf{u},t) = \left(1+\mu^{-1}\left[\left(\sum_{i=1}^{m}\lambda_i (1-u_i)+\theta\right)^{\alpha}-\theta^\alpha\right]\right)^{-\rho t},
\end{equation}	
and particularly for $\rho =1$, it satisfies
\begin{equation}\label{pgfde}
	\mu \left(1-e^{-\partial_t}\right)\mathcal{G}^{\alpha, \theta}(\textbf{u},t) = -(S(\boldsymbol{\lambda}))^\alpha \left( \left(I- \frac{\sum_{i=1}^{m}\lambda_i B_i - \theta}{S(\boldsymbol{\lambda})}\right)^\alpha -\left( \frac{\theta}{S(\boldsymbol{\lambda})}\right)^\alpha\right) \mathcal{G}^{\alpha, \theta}(\textbf{u},t).
\end{equation}
\end{theorem}

\begin{proof}
	We define the pgf of the MTSFNBP as
	\begin{equation*}
	\mathcal{G}^{\alpha, \theta}(\textbf{u},t) = \mathbb{E}[u_1^{\mathcal{Q}_1^{\alpha, \theta}(t)}\cdots u_m^{\mathcal{Q}_m^{\alpha, \theta}(t)}] = \sum_{\textbf{k} \geq \textbf{0}} u_1^{k_1} u_2^{k_2}\cdots u_m^{k_m} \mathcal{P}^{\alpha, \theta}(\textbf{k},t).
	\end{equation*}
	Using the conditioning arguments and (\ref{gammalt}), we get
	\begin{align*}
	\mathbb{E}[u_1^{\mathcal{Q}_1^{\alpha, \theta}(t)}\cdots u_m^{\mathcal{Q}_m^{\alpha, \theta}(t)}] 
	&= \mathbb{E}\left[\mathbb{E}\left[ u_1^{\mathcal{N}_1^{\alpha, \theta}(\Gamma(t),\lambda_1)}u_2^{\mathcal{N}_2^{\alpha, \theta}(\Gamma(t),\lambda_2)}\cdots u_m^{\mathcal{N}_m^{\alpha, \theta}(\Gamma(t),\lambda_m)} \;\vline\; \Gamma(t) \right]\right]\\
	&= \mathbb{E}\left[\exp\left\{\displaystyle
	-\Gamma(t)\left(\left[\sum_{i=1}^{m}\lambda_i(1-u_i) + \theta\right]^{\alpha} - \theta^\alpha\right)\right\}\right] \\
	&= \left(1+\mu^{-1}\left[\left(\sum_{i=1}^{m}\lambda_i (1-u_i)+\theta\right)^{\alpha}-\theta^\alpha\right]\right)^{-\rho t}.
	\end{align*}
Hence, we get (\ref{mtsfpgf}).\\
To prove (\ref{pgfde}), we multiply (\ref{pmf21}) by $\prod_{j=1}^{m}u_j^{k_j}$ and sum over $\textbf{k} \geq \textbf{0}$ so we get from right hand side (RHS) of (\ref{pmf21}) as
\begin{align}\label{pgfde1}
&-(S(\boldsymbol{\lambda}))^\alpha \sum_{\textbf{k} \geq \textbf{0}}\prod_{j=1}^{m}u_j^{k_j}\left( \left(I- \frac{\sum_{i=1}^{m}\lambda_i B_i - \theta}{S(\boldsymbol{\lambda})}\right)^\alpha -\left( \frac{\theta}{S(\boldsymbol{\lambda})}\right)^\alpha\right) \mathcal{P}^{\alpha, \theta}(\textbf{k},t) \nonumber \\
&= -(S(\boldsymbol{\lambda}))^\alpha \sum_{\textbf{k} \geq \textbf{0}}\prod_{j=1}^{m}u_j^{k_j}\left( \sum_{r=0}^{\infty}\binom{\alpha}{r} \left(I+\frac{\theta}{S(\boldsymbol{\lambda})}\right)^{\alpha-r} \left(\frac{-\sum_{i=1}^{m}\lambda_i B_i}{S(\boldsymbol{\lambda})}\right)^r -\left( \frac{\theta}{S(\boldsymbol{\lambda})}\right)^\alpha\right) \mathcal{P}^{\alpha, \theta}(\textbf{k},t) \nonumber\\
&= -(S(\boldsymbol{\lambda}))^\alpha \sum_{\textbf{k} \geq \textbf{0}} \sum_{r=0}^{\infty}\binom{\alpha}{r} \left(\frac{-1}{S(\boldsymbol{\lambda})}\right)^r \left(I+\frac{\theta}{S(\boldsymbol{\lambda})}\right)^{\alpha-r} \sum_{\textbf{l} \geq  \textbf{0}}\binom{r}{l_1,\ldots,l_m} \prod_{j=1}^m \lambda_j^{l_j}B_j^{l_j}u_j^{k_j}\mathcal{P}^{\alpha, \theta}(\textbf{k},t)
 -\theta^\alpha 	\mathcal{G}^{\alpha, \theta}(\textbf{u},t) \nonumber\\
 &= -(S(\boldsymbol{\lambda}))^\alpha  \sum_{r=0}^{\infty}\binom{\alpha}{r} \left(\frac{-1}{S(\boldsymbol{\lambda})}\right)^r \left(I+\frac{\theta}{S(\boldsymbol{\lambda})}\right)^{\alpha-r} \sum_{\textbf{l} \geq  \textbf{0}}\binom{r}{l_1,\ldots,l_m} \prod_{j=1}^{m}\lambda_j^{l_j}u_j^{l_j} 	\mathcal{G}^{\alpha, \theta}(\textbf{u},t)
  -\theta^\alpha 	\mathcal{G}^{\alpha, \theta}(\textbf{u},t) \nonumber\\
  &= -(S(\boldsymbol{\lambda}))^\alpha \left( \left(I- \frac{\sum_{i=1}^{m}\lambda_i u_i - \theta}{S(\boldsymbol{\lambda})}\right)^\alpha -\left( \frac{\theta}{S(\boldsymbol{\lambda})}\right)^\alpha\right) \mathcal{G}^{\alpha, \theta}(\textbf{u},t).
\end{align}
It is clear from the l.h.s. of (\ref{pmf21}) that when $\rho=1$, we have\\
$\mu \left(1-e^{-\partial_t}\right)\mathcal{G}^{\alpha, \theta}(\textbf{u},t) $
\begin{align*}
	&= \mu\mathcal{G}^{\alpha, \theta}(\textbf{u},t) - \mu\mathcal{G}^{\alpha, \theta}(\textbf{u},t-1)\\
	&= \mu\left(1+\mu^{-1}\left[\left(\sum_{i=1}^{m}\lambda_i (1-u_i)+\theta\right)^{\alpha}-\theta^\alpha\right]\right)^{-t} - \mu \left(1+\mu^{-1}\left[\left(\sum_{i=1}^{m}\lambda_i (1-u_i)+\theta\right)^{\alpha}-\theta^\alpha\right]\right)^{1-t}\\
	&= \mathcal{G}^{\alpha, \theta}(\textbf{u},t)\left[\left(\sum_{i=1}^{m}\lambda_i (1-u_i)+\theta\right)^{\alpha}-\theta^\alpha\right],
\end{align*}
which coincides with (\ref{pgfde1}). Hence, we proved (\ref{pgfde}).
\end{proof}
\noindent In the next theorem, we derive the L\'{e}vy measure density of the MTSFNBP.

\begin{theorem}
	The L\'{e}vy measure density $\Pi(\cdot)$ associated with MTSFNBP is given by
	\begin{equation}\label{lmd11}
	\Pi(\cdot) = \sum_{\textbf{k} \geq \textbf{0}}
	\left(\frac{-1}{S(\boldsymbol{\lambda})}\right)^{S(\boldsymbol{k})} \prod_{j=1}^m \frac{\lambda_j^{k_j}}{k_j !} \sum_{i=0}^{\infty} \frac{\theta^i}{i! (S(\boldsymbol{\lambda}))^i } \delta_{\{\textbf{k}\}}(\cdot) \;_2 \psi_1  \left[\frac{-(S(\boldsymbol{\lambda}))^\alpha}{\mu - \theta^{\alpha}} \; \vline \;\begin{matrix}
	\left(1,\alpha\right), (0, 1)\\
	(1-S(\boldsymbol{k})-i, \alpha)
	\end{matrix} \right],
	\end{equation}
	which is absolutely convergent for $ \mu > (S(\boldsymbol{\lambda}))^\alpha +\theta^ \alpha.$ Here $\delta_{\{\textbf{k}\}}(\cdot)$  is the Dirac measure concentrated at $\textbf{k}$. 
\end{theorem}
\begin{proof}
	Let $\mu(s) = \rho s^{-1}e^{-\mu s}$ is the L\'{e}vy density of the gamma subordinator. Then,
	applying Theorem 30.1 in \cite{Sato} and with the help of (\ref{pmf2}), we get
	\begin{align*}
		\Pi(\cdot) 
		&= \int_{0}^{\infty} \sum_{\textbf{k} \geq \textbf{0}} 	q^{\alpha, \theta}(\textbf{k},s) \delta_{\{k\}}(\cdot) \mu(s) ds\\
		&= 	 \sum_{\textbf{k} \geq \textbf{0}} \left(\frac{-1}{S(\boldsymbol{\lambda})}\right)^{S(\boldsymbol{k})} \prod_{j=1}^m \frac{\lambda_j^{k_j}}{k_j !} \sum_{i=0}^{\infty} \rho \theta^i  \delta_{\{\textbf{k}\}}(\cdot)\sum_{r \geq 0}\binom{\alpha r}{i} \binom{\alpha r-i}{S(\textbf{k})}\frac{(S(\boldsymbol{\lambda}))^{\alpha r-i}}{r!}\int_{0}^{\infty} s^{r-1} e^{s(\theta^\alpha - \mu)}ds\\
			&= 	 \sum_{\textbf{k} \geq \textbf{0}} \left(\frac{-1}{S(\boldsymbol{\lambda})}\right)^{S(\boldsymbol{k})} \prod_{j=1}^m \frac{\lambda_j^{k_j}}{k_j !} \sum_{i=0}^{\infty} \rho \theta^i \delta_{\{\textbf{k}\}}(\cdot) \sum_{r \geq 0}\binom{\alpha r}{i} \binom{\alpha r-i}{S(\textbf{k})}\frac{(S(\boldsymbol{\lambda}))^{\alpha r-i}}{r!} \frac{\Gamma(r)}{(\mu-\theta^\alpha)^r}.
	\end{align*}
	With the help of the generalized Wright function, we obtain the  L\'{e}vy measure density $\Pi$.
\end{proof}

\begin{remark}
	When $\theta =0$, (\ref{lmd11}) reduces to 
	\begin{equation*}
		\Pi(\cdot) = \sum_{\textbf{k} \geq \textbf{0}}
	\left(\frac{-1}{S(\boldsymbol{\lambda})}\right)^{S(\boldsymbol{k})} \prod_{j=1}^m \frac{\lambda_j^{k_j}}{k_j !} \delta_{\{\textbf{k}\}}(\cdot) \;_2 \psi_1  \left[\frac{-(S(\boldsymbol{\lambda}))^\alpha}{\mu} \; \vline \;\begin{matrix}
	\left(1,\alpha\right), (0, 1)\\
	(1-S(\boldsymbol{k}), \alpha)
	\end{matrix} \right],
	\end{equation*}
	which is same as L\'{e}vy measure density of the multivariate space-fractional negative binomial process studied in \cite{Beghin2018}.
\end{remark}

\section{Bivariate Risk Models with shocks}
Let $\Gamma(t) \sim G(\mu, \rho t)$ be the gamma subordinator. Let $\mathcal{Q}_0(t), \mathcal{Q}_1(t)$ and $\mathcal{Q}_2(t)$ are three i.i.d. TSFNBPs such that $\mathcal{Q}_i(t) = \mathcal{N}_i^{\alpha, \theta}(\Gamma(t), \lambda_i)$ with $\lambda_i >0,  \;\;  i=0,1,2$.\\
 We define the counting process as
\begin{equation}\label{bcp1}
N_1(t) = \mathcal{Q}_1(t) + \mathcal{Q}_0(t), \;\;\; N_2(t) = \mathcal{Q}_2(t) + \mathcal{Q}_0(t).
\end{equation}
Let $\nu \geq 0$ be the initial capital and $\eta > 0$ be the premium rate per unit time. Consider the surplus risk model based on TSFNBP defined as
\begin{equation}\label{brm1}
	R_{\mathcal{Q}}(t) = \nu +\omega t - S(t),\;\;\; t \geq 0,
\end{equation}
where $\{S(t)\}_{t \geq 0}$ be the total claim amount process defined by
\begin{equation}
S(t) = S_1(t) + S_2(t),
\end{equation}
with
\begin{align*}
S_1(t) &= \sum_{i=0}^{\mathcal{Q}_1(t)} \xi_{1,i} + \sum_{i=0}^{\mathcal{Q}_0(t)} \xi_{3,i}, \;\; \text{with}\;\; \xi_{1,0} = \xi_{3,0} =0, \\
S_2(t) &= \sum_{i=0}^{\mathcal{Q}_2(t)} \xi_{2,i} + \sum_{i=0}^{\mathcal{Q}_0(t)} \xi_{4,i}, \;\; \text{with}\;\; \xi_{2,0} = \xi_{4,0} =0.
\end{align*}

The sequence of random variables $\xi_{k,i}, \; k=1,2,3,4$ and $\mathcal{Q}_j(t),\; j=0,1,2$ are mutually independent. We may observe that the counting processes $N_1(t)$ and $N_2(t)$ are dependent. Consequently, $S_1(t)$ and $S_2(t)$ are generally dependent.

\begin{remark}
	This bivariate risk model allows the claim to be of two mutually exclusive different types and $\mathcal{Q}_0(t)$ counts the numbers of shock events. The distribution of $\xi_{k,i}$ describes the $i$th claim amount at certain type $k$.
\end{remark}

\subsection{The Bivariate Counting Process} \hfill \\
Let $\{(N_1(t), N_2(t)):  t \geq 0\}$ be the bivariate counting process (BCP) based on the  TSFNBP as defined in (\ref{bcp1}). With the help of (\ref{mtsfpgf}) for $m=3$, we have
\begin{align}\label{bpgf1}
\mathbb{E}[u_1^{N_1(t)}u_2^{N_2(t)}] 
&= \mathbb{E}[u_1^{\mathcal{Q}_1(t) + \mathcal{Q}_0(t)}u_2^{\mathcal{Q}_2(t) + \mathcal{Q}_0(t)}]\nonumber\\
&= \mathbb{E}[u_1^{\mathcal{Q}_1(t)}u_2^{\mathcal{Q}_2(t)}(u_1 u_2)^{\mathcal{Q}_0(t)}]\nonumber \\
&= \sum_{k_1, k_2, k_3 \geq 0} u_1^{k_1}u_2^{k_2}(u_1 u_2)^{k_3} Pr\{\mathcal{Q}_1(t) = k_1, \mathcal{Q}_2(t) =k_2, \mathcal{Q}_0(t) =k_3\}\nonumber\\
&= \left(1+\mu^{-1}\left[\left(\lambda_1 (1-u_1)+\lambda_2 (1-u_2) + \lambda_0(1-u_1 u_2) + \theta\right)^{\alpha}-\theta^\alpha\right]\right)^{-\rho t}.
\end{align}
\noindent To fomulate following results, we first define a random vector $(b_1, b_2, b_0)$  with pmf given by
\begin{equation}\label{pmfve1}
\begin{cases}
Pr\{b_1 =1, b_2 = 0, b_0=0\} &= \frac{\lambda_1}{\lambda},\\
Pr\{b_1 =0, b_2 = 1, b_0=0\} &= \frac{\lambda_2}{\lambda},\\
Pr\{b_1 =0, b_2 = 0, b_0=1\} &= \frac{\lambda_0}{\lambda},
\end{cases}
\end{equation}
where $\lambda =\lambda_1 +\lambda_2+\lambda_0$.\\
Let $\mathcal{Q}(t) = \mathcal{N}^{\alpha, \theta}(\Gamma(t), \lambda)$ be the TSFNBP with $\mathcal{Q}(0) =0$. Let $(b_{1,i}, b_{2,i}, b_{0,i}), \; i=1,2,\ldots$ be the i.i.d random vectors independent from $\mathcal{Q}$ with pmf given in (\ref{pmfve1}) . Then, we define
\begin{equation*}
H_1(t) = \sum_{i=0}^{\mathcal{Q}(t)}(b_{1,i}+b_{0,i}), \;\;\; H_2(t) = \sum_{i=0}^{\mathcal{Q}(t)}(b_{2,i}+b_{0,i}),
\end{equation*}
with $b_{1,0} = b_{2,0} = b_{0,0} =0$.

In the next proposition, we show that the BCP is a particular case of the bivariate compound Poisson process associated with TSFNBP.

\begin{proposition}
	The BCP $(N_1, N_2)$ and $(H_1, H_2)$ are identical in distribution.
\end{proposition}

\begin{proof}
	We derive the pgf of the BCP, $(H_1, H_2)$ as
	\begin{align*}
	\mathbb{E}[u_1^{H_1(t)}u_2^{H_2(t)}] 
	&=\mathbb{E}[u_1^{\sum_{i=0}^{\mathcal{Q}(t)}(b_{1,i}+b_{0,i})}u_2^{\sum_{i=0}^{\mathcal{Q}(t)}(b_{2,i}+b_{0,i})}]\\
	&= \sum_{k=0}^{\infty} \left[\mathbb{E}[u_1^{b_{1,i}+b_{0,i}}u_2^{b_{2,i}+b_{0,i}}]\right]^k Pr\{\mathcal{Q}(t) =k\}\\
		&= \sum_{k=0}^{\infty} \left[\frac{\lambda_1}{\lambda}u_1 +\frac{\lambda_2}{\lambda}u_2+ \frac{\lambda_0}{\lambda}u_1 u_2 \right]^k Pr\{\mathcal{Q}(t) =k\}.
	\end{align*}
	With the help of (\ref{mtsfpgf}) for $m=1$, we get
	\begin{align*}
		\mathbb{E}[u_1^{H_1(t)}u_2^{H_2(t)}] 
		&= \left(1+\mu^{-1}\left[\left(\lambda \left(1-\frac{\lambda_1}{\lambda}u_1 -\frac{\lambda_2}{\lambda}u_2- \frac{\lambda_0}{\lambda}u_1 u_2\right)+\theta\right)^{\alpha}-\theta^\alpha\right]\right)^{-\rho t}\\
&= \left(1+\mu^{-1}\left[\left(\lambda_1 (1-u_1)+\lambda_2 (1-u_2) + \lambda_0(1-u_1 u_2) + \theta\right)^{\alpha}-\theta^\alpha\right]\right)^{-\rho t},
	\end{align*}
 	which coincides with (\ref{bpgf1}). Hence, $(N_1, N_2)$ and $(H_1, H_2)$ are identical in distribution.
\end{proof}
\noindent The mean and variance of the counting processes $N_i(t), \; i=1,2$ are gievn by
\begin{equation}
\mathbb{E}[N_i(t)] = \frac{\alpha \rho t \theta^{\alpha-1}}{\mu} (\lambda_i +\lambda_0),
\end{equation}
and 
\begin{equation}
\text{Var}[N_i(t)] = \left[(\lambda_i^2 +\lambda_0^2)\alpha \theta^{\alpha-2}\left(\frac{\alpha \theta^\alpha}{\mu}+1-\alpha\right) + \alpha \theta^{\alpha-1} (\lambda_i +\lambda_0)\right]\frac{\rho t}{\mu}.
\end{equation}
\begin{remark}
	It is observed the $\text{Var}[N_i(t)] - \mathbb{E}[N_i(t)] >0$. Hence, the model can accomodate the overdispersed situations.
\end{remark}
\subsection{Stochastically Equivalent Claim Amount Process} \hfill \\
Let $S_3(t)$ and $S_4(t)$ be two compound processes defined by 
\begin{equation}\label{cap1}
S_3(t) = \sum_{i=0}^{\mathcal{Q}(t)}(b_{1,i} \xi_{1,i}+b_{0,i} \xi_{3,i}), \;\;\; S_4(t) = \sum_{i=0}^{\mathcal{Q}(t)}(b_{2,i} \xi_{2,i}+b_{0,i} \xi_{4,i}),
\end{equation}
where $\mathcal{Q}(t)$ is TSFNBP with parameter $\lambda = \lambda_1 +\lambda_2+\lambda_0$. For each fix $i$, $\mathcal{Q}(t), \xi_{1,i}, (\xi_{3,i}, \xi_{4,i}), \xi_{2,i}$ and the random vector $(b_{1,i}, b_{2,i}, b_{0,i})$ are mutually independent.\\
Then, we have
\begin{align}\label{bpgf111}
	\mathbb{E}[u_1^{S_3(t)}u_2^{S_4(t)}]&=\mathbb{E}\left[u_1^{\sum_{i=0}^{\mathcal{Q}(t)}(b_{1,i} \xi_{1,i}+b_{0,i} \xi_{3,i})}u_2^{\sum_{i=0}^{\mathcal{Q}(t)}(b_{2,i} \xi_{2,i}+b_{0,i} \xi_{4,i})} \;\vline\; \mathcal{Q}(t) \right] \nonumber \\
&= \sum_{n=0}^{\infty} \mathbb{E}\left[u_1^{\sum_{i=0}^{n}(b_{1,i} \xi_{1,i}+b_{0,i} \xi_{3,i})}u_2^{\sum_{i=0}^{n}(b_{2,i} \xi_{2,i}+b_{0,i} \xi_{4,i})}\right] Pr\{\mathcal{Q}(t) =n\}\nonumber \\
&= \sum_{n=0}^{\infty} \left(\mathbb{E}\left[\mathbb{E}\left[u_1^{b_{1,1} \xi_{1,i}}u_2^{b_{2,1} \xi_{2,1}} (u_1^{\xi_{3,1}} u_2^{\xi_{4,1}})^{b_{0,1}} \; \vline \; ((b_{1,1}, b_{2,1}, b_{0,1}))\right]\right]\right)^n Pr\{\mathcal{Q}(t) =n\}\nonumber \\
&= \sum_{n=0}^{\infty} \left( \frac{\lambda_1}{\lambda} \mathbb{E} u_1^{\xi_{1,1}} + \frac{\lambda_2}{\lambda} \mathbb{E}u_2^{\xi_{2,1}} +\frac{\lambda_0}{\lambda} \mathbb{E}(u_1^{\xi_{3,1}} u_2^{\xi_{4,1}})\right)^n Pr\{\mathcal{Q}(t) =n\}\nonumber \\
 &=  \left(1+\mu^{-1}\left[\left(\lambda_1 (1-\mathbb{E} u_1^{\xi_{1,1}})+\lambda_2 (1- \mathbb{E}u_2^{\xi_{2,1}}) + \lambda_0(1-\mathbb{E}(u_1^{\xi_{3,1}} u_2^{\xi_{4,1}}) + \theta\right)^{\alpha}-\theta^\alpha\right]\right)^{-\rho t}.
\end{align}
Using the above construction, we prove the following proposition.
\begin{proposition}
	The claim amount process $(S_1, S_2)$ is stochastically equivalent to the process $(S_3, S_4)$ defined in (\ref{cap1}).
\end{proposition}
\begin{proof}
	We compute the pgf of the process $(S_1, S_2)$ as\\
	$		\mathbb{E}[u_1^{S_1(t)}u_2^{S_2(t)}]  $
	\begin{align*}
		&=\mathbb{E}\left[u_1^{\sum_{i=0}^{\mathcal{Q}_1(t)} \xi_{1,i} + \sum_{i=0}^{\mathcal{Q}_0(t)} \xi_{3,i}}u_2^{\sum_{i=0}^{\mathcal{Q}_2(t)} \xi_{2,i} + \sum_{i=0}^{\mathcal{Q}_0(t)} \xi_{4,i}}\right]\\
		&= \sum_{k_1, k_2, k_3 \geq 0} \mathbb{E}u_1^{\sum_{i=0}^{\mathcal{Q}_1(t)} \xi_{1,i}} \mathbb{E}u_2^{\sum_{i=0}^{\mathcal{Q}_2(t)} \xi_{2,i}} \mathbb{E}[u_1^{\sum_{i=0}^{\mathcal{Q}_0(t)} \xi_{3,i}}u_2^{\sum_{i=0}^{\mathcal{Q}_0(t)} \xi_{4,i}}]Pr\{\mathcal{Q}_1(t) = k_1, \mathcal{Q}_2(t) =k_2, \mathcal{Q}_0(t) =k_3\}\\
			&= \sum_{k_1, k_2, k_3 \geq 0} \left(\mathbb{E}u_1^{\xi_{1,i}}\right)^{k_1} \left(\mathbb{E}u_2^{ \xi_{2,i}}\right)^{k_2} \left(\mathbb{E}[u_1^{ \xi_{3,i}}u_2^{\xi_{4,i}}]\right)^{k_3}Pr\{\mathcal{Q}_1(t) = k_1, \mathcal{Q}_2(t) =k_2, \mathcal{Q}_0(t) =k_3\}.
	\end{align*}
	 Using (\ref{mtsfpgf}) for $m=3$, we get\\
	 	$\mathbb{E}[u_1^{S_1(t)}u_2^{S_2(t)}]$
	 \begin{equation*}
	 =  \left(1+\mu^{-1}\left[\left(\lambda_1 (1-\mathbb{E} u_1^{\xi_{1,1}})+\lambda_2 (1- \mathbb{E}u_2^{\xi_{2,1}}) + \lambda_0(1-\mathbb{E}(u_1^{\xi_{3,1}} u_2^{\xi_{4,1}}) + \theta\right)^{\alpha}-\theta^\alpha\right]\right)^{-\rho t},
	 \end{equation*}
	 which coincides with (\ref{bpgf111}). Hence, the proposition is proved.
\end{proof}
\noindent Consequently, we state the following corollary.
\begin{corollary}\label{c1}
	The pgf of the time intersection 
	\begin{equation}
	S(t) = \sum_{i=0}^{\mathcal{Q}(t)} b_{1,i} \xi_{1,i}+b_{0,i} (\xi_{3,i} + \xi_{4,i}) + b_{2,i} \xi_{2,i}, \;\; t \geq 0,
	\end{equation}
	of the compound Poisson process $S(t)$ takes the following form
	\begin{equation*}
	\mathbb{E}u^{S(t)}=  \left(1+\mu^{-1}\left[\left(\lambda_1 (1-\mathbb{E} u^{\xi_{1,1}})+\lambda_2 (1- \mathbb{E}u^{\xi_{2,1}}) + \lambda_0(1-\mathbb{E}u^{{\xi_{3,1}} +{\xi_{4,1}}} + \theta\right)^{\alpha}-\theta^\alpha\right]\right)^{-\rho t}.
	\end{equation*} 
\end{corollary}

Now, we propose  an univariate risk model based on the TSFNBP  which is stochastically equivalent to the bivariate risk model in (\ref{brm1}).
\subsection{Stochastically Equivalent Risk Models} \hfill \\
Let $\mathcal{Q}(t) = \mathcal{N}(S^{\alpha, \theta}(\Gamma(t)), \lambda)$ be a TSFNBP with parameter $\lambda = \lambda_1 +\lambda_2+\lambda_0$. Then, we consider the surplus risk model governed by the TSFNBP defined as
\begin{equation}\label{sm1}
R_{\mathcal{Q}}(t) = \nu + \omega t- \sum_{i=0}^{\mathcal{Q}(t)}\varphi_i, \;\; t \geq 0,
\end{equation}
where $\varphi_i =b_{1,i} \xi_{1,i}+b_{0,i} (\xi_{3,i} + \xi_{4,i}) + b_{2,i} \xi_{2,i}, \;\varphi_0=0  $. Here $\omega >0$ denotes the constant premium rate and $\nu>0$ is the initial capital.\\
Let  $\varphi_i, i=1,2,\ldots$ be i.i.d. claim amounts with commom distribution function $F_\varphi$. Then, the premium loading factor $\theta$ which presents the insurance firm's profit margin is given by

\begin{equation}
\theta = \frac{\mathbb{E}\left[ \omega t- \sum_{i=0}^{\mathcal{Q}(t)}\varphi_i\right]}{ \mathbb{E}\left[\sum_{i=0}^{\mathcal{Q}(t)}\varphi_i\right]} = \frac{\omega t}{\mathbb{E}[\mathcal{Q}(t) ] \mathbb{E}[\varphi_i]} -1 = \frac{\omega t}{\alpha  \lambda \rho t \theta^{\alpha -1}\mu^{-1}\mathbb{E}[\varphi_i]} -1,
\end{equation} 
where $ \mathbb{E}[\varphi_1] = \frac{\lambda_1}{\lambda} \mathbb{E} [\xi_{1,1}] + \frac{\lambda_2}{\lambda} \mathbb{E}[\xi_{2,1}] +\frac{\lambda_0}{\lambda} \mathbb{E}[\xi_{3,1} + \xi_{4,1}].$\\
The net profit condition is 
\begin{equation*}
\theta > 0 \;\; \Longrightarrow \;\; \omega > \lambda_1 \mathbb{E} [\xi_{1,1}] + \lambda_2 \mathbb{E}[\xi_{2,1}] +\lambda_0\mathbb{E}[\xi_{3,1} + \xi_{4,1}].
\end{equation*}
 Let $T_R$ denotes the first-hitting time to ruin defined by
\begin{equation}
T_R = \inf \{t > 0 : R_{\mathcal{Q}}(t) <0\},
\end{equation}
with probability $\mathcal{P}(u) = Pr\{T_R < \infty \}$. The joint probability $J(u,y)$ that the ruin occurs at finite time and the deficit at the time of ruin denoted by $\mathbb{D} = |R_{\mathcal{Q}}(t)|$ is defined as
\begin{equation}
J(u,y) = Pr\{T_R < \infty, \mathbb{D} \leq y\}, \;\; y \geq 0.
\end{equation}
First, we state the following lemma to prove the subsequent results.
{\begin{lemma}\label{lemma11}
		The transition probabilities of the TSFNBP are given by (see \cite{Orsingher2015})
		\begin{equation}
		Pr\{\mathcal{Q}(t+h) =  n | \mathcal{Q}(t) =m\}\nonumber =\left\{
		\begin{array}{ll}
		1-h \psi(\lambda) +o(h) & \quad n = m, \\
		-h\frac{(-\lambda)^j}{j!}\psi^{(j)}(\lambda) + o(h)  & \quad n = m+j,\; j=1,2,3,\dots,
		\end{array}
		\right.
		\end{equation}
		where $\psi(\lambda) = \rho \left[\log \left(\mu -\theta^{\alpha} + (\theta + \lambda)^\alpha\right) - \log (\alpha)\right] $ is the Laplace exponent of $ \{S^{\alpha, \theta}(\Gamma(t))\}_{t \geq 0}$ (see \cite[Equation (3.3)]{Kumar2019}).
\end{lemma}
Now, we discuss the ruin probability, joint distribution of time to ruin and deficit at ruin and the governing integro-differential equation under the assumptions of the model in (\ref{sm1}).
\begin{theorem}
	The joint distribution $J(u,y)$ satisfies the following integro-differential equation.
	\begin{equation}\label{ju1}
	\diffp{J}{u} = \frac{\psi(\lambda)}{\omega} \left(J(u, y) - \int_{0}^{u}J(u-x, y)dF_{\varphi} + F_\varphi(u+y)-F_\varphi(u)\right).
	\end{equation}
\end{theorem}
\begin{proof}
	Using Lemma \ref{lemma11}, we get
	\begin{align*}
	J(u,y) = (1&-h\psi(\lambda)) J(\hat{u},y) - h\sum_{j = 1}^{\infty} \frac{(-\lambda)^j}{j!}\psi^{(j)}(\lambda)\\
	& \times \left(\int_{0}^{\hat{u}}J(\hat{u}-x, y)dF_{\varphi} + F_\varphi(\hat{u}+y)-F_\varphi(\hat{u})\right),
	\end{align*}
	where $\hat{u} = u+\omega t$.
	\\ On rearranging terms, we get
	\begin{align*}
	\frac{J(\hat{u},y) - J(u,y)}{\omega h} &= \frac{\psi(\lambda)}{\omega} J(\hat{u}, y) + \frac{1}{\omega} \sum_{j = 1}^{\infty} \frac{(-\lambda)^j}{j!}\psi^{(j)}(\lambda)\\ & \times \left(\int_{0}^{\hat{u}}J(\hat{u}-x, y)dF_{\varphi} + F_\varphi(\hat{u}+y)-F_\varphi(\hat{u})\right).
	\end{align*}
	Using Taylor series and limiting $h \rightarrow 0$, we get
	\begin{equation}
	\diffp{J}{u} = \frac{\psi(\lambda)}{\omega} J(u, y) - \frac{\psi(\lambda)}{\omega} \left(\int_{0}^{u}J(u-x, y)dF_{\varphi} + F_\varphi(u+y)-F_\varphi(u)\right).
	\end{equation}
	On simplifying, we get the required integro-differential equation.
\end{proof}
\begin{remark}
	For $y \rightarrow \infty$, we get 
	\begin{equation}
	\frac{d \mathcal{P}(u)}{du} = \frac{\psi(\lambda)}{\omega} \left( \mathcal{P}(u)+ 1 - F_{\varphi}(u) + \int_{0}^{u}\mathcal{P}(u-x)dF_{\varphi} \right).
	\end{equation}
\end{remark}
\begin{theorem}
	The joint distribution $J(u,y)$ when the initial capital is zero, is given by
	\begin{equation*}
	J(0,y) = \frac{\psi(\lambda)}{\omega} \int_{0}^{\infty} ( F_\varphi(u+y)-F_\varphi(u))du.
	\end{equation*}
\end{theorem}
\begin{proof}
	Integrating both sides of (\ref{ju1}), we get
	\begin{equation*}
	J(0,y) = J(0, \infty) + \frac{\psi(\lambda)}{\omega} \int_{0}^{\infty} \left(J(u, y) - \int_{0}^{u}J(u-x, y)dF_{\varphi} + F_\varphi(u+y)-F_\varphi(u)\right)du.
	\end{equation*}
	Using $J(\infty ,y) =0$, we get the result. 
\end{proof}
\begin{remark}
	It is noted that when $y \rightarrow \infty$, we get
	\begin{equation}
	\mathcal{P}(0) = \frac{\psi(\lambda)}{\omega} \int_{0}^{\infty} (1-F_\varphi(u)) du.
	\end{equation}
\end{remark}

\subsubsection{\textbf{Long-Range Dependence of the Risk Process in (\ref{sm1})}}

\begin{proposition}
	Under the assumptions of the model in (\ref{sm1}), the covariance structure of $R_{\mathcal{Q}}(t)$ is given by
	\begin{equation*}
	\text{Cov}[R_{\mathcal{Q}}(t), R_{\mathcal{Q}}(s)] =  \alpha \rho s \mu^{-1} \theta^{\alpha-1} \text{Var}[\varphi_1] + (\mathbb{E}[\varphi_1])^2 (\alpha (1-\alpha) \rho \mu^{-1} \theta^{\alpha-2} + \rho \mu^{-2}\alpha^{2} \theta^{2\alpha -2})s.
	\end{equation*}
	Also, the variance of $R_{\mathcal{Q}}(t)$ is 
		\begin{equation*}
	\text{Var}[R_{\mathcal{Q}}(t)] = \alpha \rho t \mu^{-1} \theta^{\alpha-1} \text{Var}[\varphi_1] + (\mathbb{E}[\varphi_1])^2 (\alpha (1-\alpha) \rho \mu^{-1} \theta^{\alpha-2} + \rho \mu^{-2}\alpha^{2} \theta^{2\alpha -2})t.
	\end{equation*}
\end{proposition}
\begin{proof}
	For $0 < s \leq t$, the covariance of $R_{\mathcal{Q}}(t)$ can be obtained along similar lines to that of Proposition 10 of \cite{Kumar2020} as
	\begin{align*}
		\text{Cov}[R_{\mathcal{Q}}(t), R_{\mathcal{Q}}(s)] 
		&= \text{Cov}\left[\sum_{i=0}^{\mathcal{Q}(t)}\varphi_i, \sum_{i=0}^{\mathcal{Q}(s)}\varphi_i\right]\\
		&= \text{Var}[\varphi_1] \mathbb{E}[S^{\alpha, \theta}(\Gamma(t))] + (\mathbb{E}[\varphi_1])^2 \text{Cov}[S^{\alpha, \theta}(\Gamma(t)), S^{\alpha, \theta}(\Gamma(s))]\\
		&= \alpha \rho s \mu^{-1} \theta^{\alpha-1} \text{Var}[\varphi_1] + (\mathbb{E}[\varphi_1])^2 (\alpha (1-\alpha) \rho \mu^{-1} \theta^{\alpha-2} + \rho \mu^{-2}\alpha^{2} \theta^{2\alpha -2})s.
	\end{align*}
	The last line is due to Section 3.3 and Section 3.5 of \cite{Kumar2019}.\\
	When $s=t$, $	\text{Cov}[R_{\mathcal{Q}}(t), R_{\mathcal{Q}}(s)] = \text{Var}[R_{\mathcal{Q}}(t)]$. Hence, the proposition follows.
\end{proof}
\noindent We now prove the LRD property of the $R_{\mathcal{Q}}(t)$.
\begin{proposition}
	The risk process $R_{\mathcal{Q}}(t)$ has the LRD property for $\alpha \in (0,1)$.
\end{proposition}

\begin{proof}
		The correlation function of the $R_{\mathcal{Q}}(t)$ is
	\begin{align*}
	\text{Corr}[R_{\mathcal{Q}}(t) , R_{\mathcal{Q}}(s) ] &= \frac{\text{Cov}[R_{\mathcal{Q}}(t) , R_{\mathcal{Q}}(s)]}{\sqrt{\text{Var}[R_{\mathcal{Q}}(t)]}\sqrt{\text{Var} [R_{\mathcal{Q}}(s)]}}\\
	&=  \frac{ \alpha \rho s \mu^{-1} \theta^{\alpha-1} \text{Var}[\varphi_1] + (\mathbb{E}[\varphi_1])^2 (\alpha (1-\alpha) \rho \mu^{-1} \theta^{\alpha-2} + \rho \mu^{-2}\alpha^{2} \theta^{2\alpha -2})s}{\sqrt{st}\sqrt{ [\alpha \rho  \mu^{-1} \theta^{\alpha-1} \text{Var}[\varphi_1] + (\mathbb{E}[\varphi_1])^2 (\alpha (1-\alpha) \rho \mu^{-1} \theta^{\alpha-2} + \rho \mu^{-2}\alpha^{2} \theta^{2\alpha -2})]^2}}\\
	&=s^{1/2}t^{-1/2}.
	\end{align*}
	Hence, we get
	\begin{equation*}
	\lim_{t \rightarrow \infty} \frac{\text{Corr}[R_{\mathcal{Q}}(t) , R_{\mathcal{Q}}(s)]}{t^{-1/2}} = s^{1/2},
	\end{equation*}
	which shows that $d=\frac{1}{2} \in (0,1)$. Therefore, $R_{\mathcal{Q}}(t)$ exhibits the LRD property.
\end{proof}

\section{Concluding Remarks}
In this study, we have defined a multivariate tempered space fractional Poisson process and discussed its important properties. We explored the mixture tempered variants of the MTSFPP and their connection to PDEs. Using an independent gamma subordinator, we proposed a multivariate tempered space fractional negative binomial process  (MTSFNBP) and studied its distributional properties along with its L\'{e}vy measure density.  The MTSFNBP generalizes the TSFNBP and the multivariate space fractional negative binomial process investigated by Maheshwari \cite{Maheshwari2023} and Beghin and Vellaisamy \cite{Beghin2018}, respectively. We presented a bivariate risk model based on the TSFNBP and obtained its stochastically equivalent univariate generalized Cramer-Lundberg risk process.  We obtained several ruin-related measures and established the LRD property for the risk process. In insurance, our approach can be useful when multiple types of insurance coverage and shock events affect the total claim amount.



		\end{document}